\documentclass[letterpaper, 10 pt, conference]{ieeeconf}  % Comment this line out if you need a4paper

\IEEEoverridecommandlockouts                              % This command is only needed if 
\title{\LARGE \bf
A Distributed Linear Quadratic Discrete-Time Game Approach\\  to Formation Control with Collision Avoidance 
}

\author{Prima Aditya and Herbert Werner% <-this % stops a space
\thanks{Prima Aditya and Herbert Werner are with Institute of Control Systems, 
        Hamburg University of Technology, 21073, Hamburg, Germany
        {\tt\small \{prima.aditya, h.werner\}@tuhh.de}}%
}

\usepackage{fancyhdr,lipsum}
\makeatletter
\fancypagestyle{firstpage}{% Page style for first page
	  \fancyhf{}% Clear header/footer
	  \fancyhead[C]{This document is the author's version of the full paper, and should the paper be published, it will be copyrighted by IEEE.}% Header
	}
\begin{document}

\maketitle
%% Include the next line when compiling for arxiv upload
\thispagestyle{firstpage}% firstpage page style for first page

%%%%%%%%%%%%%%%%%%%%%%%%%%%%%%%%%%%%%%%%%%%%%%%%%%%%%%%%%%%%%%%%%%%%%%%%%%%%%%%%
\begin{abstract}
Formation control problems can be expressed as linear quadratic discrete-time games (LQDTG) for which Nash equilibrium solutions are sought. However, solving such problems requires solving coupled Riccati equations, which cannot be done in a distributed manner. A recent study showed that a distributed implementation is possible for a consensus problem when fictitious agents are associated with edges in the network graph rather than nodes. This paper proposes an extension of this approach to formation control with collision avoidance, where collision is precluded by including appropriate penalty terms on the edges. To address the problem, a state-dependent Riccati equation needs to be solved since the collision avoidance term in the cost function leads to a state-dependent weight matrix.  This solution provides relative control inputs associated with the edges of the network graph. These relative inputs then need to be mapped to the physical control inputs applied at the nodes;  this can be done in a distributed manner by iterating over a gradient descent search between neighbors in each sampling interval.  Unlike inter-sample iteration frequently used in distributed MPC, only a matrix-vector multiplication is needed for each iteration step here, instead of  an optimization problem to be solved.  This approach can be implemented in a receding horizon manner, this is demonstrated through a numerical example.
\end{abstract}

%%%%%%%%%%%%%%%%%%%%%%%%%%%%%%%%%%%%%%%%%%%%%%%%%%%%%%%%%%%%%%%%%%%%%%%%%%%%%%%%
\section{INTRODUCTION}

Distributed control of multi-agent (in the sense of multi-vehicle) systems has been extensively studied over the last two decades with potential applications in many areas. Formation control is one such problem that has received significant attention. In formation control, all agents in a multi-agent system must move from  arbitrary initial states to attain a pre-determined geometric shape \cite{Chen}. To attain and maintain the formation, the agents in the team exchange information about their positions and velocities.

% Numerous research topics within the multi-agent system often study formation control that is modeled by combining individual agent objectives to create a global formation objective. However, 

When formation control schemes are implemented in a distributed manner, then in situations that involve e.g.\  collision avoidance, agents may have conflicting interests, and achieving their individual objectives may take precedence over cooperation. Such situations reflect non-cooperative game behavior, as  agents strive to meet their goals without collaboration. The solution to this type of game is to find a Nash equilibrium, where individual agents cannot improve their payoff by changing their strategy unilaterally. Linear quadratic differential games (LQDG) have been proposed as a means of addressing this problem, where the cost of each agent is quadratic, and agent dynamics are assumed to be linear.

A formation control problem modeled as LQDG has been discussed in \cite{Gu}. There, a coupled Riccati differential equation is solved to find a Nash equilibrium. A discrete-time version of LQDG, referred to as  linear quadratic discrete-time game (LQDTG), is more appropriate for receding horizon implementations, and has been proposed in \cite{Jond1}. However, solving coupled Riccati differential or difference equations is likely to be intractable for large networks; moreover, the solution cannot be implemented in a distributed manner. 

Based on an idea proposed in \cite{Jond2}, it was shown in \cite{Prima} that one can avoid solving coupled Riccati difference equations by relocating the coupling terms that initially appear in the cost function to the system dynamics. Consequently, the modified problem can be reformulated as a fictitious multi-agent system evolving on the edges of the network graph instead of the nodes, allowing for a distributed solution to the decoupled Riccati difference equations. The resulting relative control inputs associated with each edge can then be mapped back to the physical control inputs in a distributed manner by employing a distributed steepest descent iteration between agents over two sampling instants. Such intersample iteration is frequently used in distributed MPC \cite{Stomberg}. However, unlike in distributed MPC, the approach proposed in \cite{Prima} does not require solving an optimization problem at each iteration step but only involves performing a matrix-vector multiplication.

Whereas the distributed scheme proposed in \cite{Prima} considers an unconstrained consensus control problem, our contribution in this article is to extend this approach to a formation control problem that includes collision avoidance among agents. We begin by formulating the problem on the graph nodes, considering the desired formation displacements together with relative constraints for collision avoidance which are represented as soft constraints in the cost function. These state-dependent collision avoidance terms in the cost lead to coupled state-dependent Riccati difference equations (SDRDE). To decouple these, we use the same idea as in \cite{Prima} by relocating the coupling term from the cost function to the system dynamics on the edges of the graph. This results in a decoupled cost that still incorporates the collision avoidance term. The reformulated problem involves solving a set of decoupled SDRDEs. This can be achieved using a receding horizon technique proposed in \cite{Dutka}, and can be implemented in a distributed manner.

% The cost of each fictitious edge agent is independent of the others, making it possible to consider a global objective by adding all individual edge costs. Instead of a Nash equilibrium, the edge problem now reflects an optimization problem with a Pareto optimal solution \cite{Reddy}. In this article, we  demonstrate a comparison of the global cost incurred when solving the coupled problem, with the solution obtained by  running the problem on the graph's edges.

The paper is organized as follows: Section II provides a review of graph theory and the formation control problem with collision avoidance.  Our proposed distributed solution is outlined in Section III. Section IV showcases simulation results, and finally, Section V concludes this article.

\section{PRELIMINARIES}

\subsection{Graph Theory}
A  graph $\mathcal{G} := (\mathcal{V},\mathcal{E})$ consists of a set of nodes $\mathcal{V}=\{\nu^1,...,\nu^N\}$, and a set of edges $\mathcal{E}=\{(\nu^i,\nu^j)\in \mathcal{V}\times \mathcal{V}, \nu^j\neq \nu^i\}$ which contains  ordered pairs of distinct nodes. $N$ is the number of nodes, and $M$ is the number of edges.  $\mathcal{G}$ is called undirected if $(\nu^i,\nu^j)\in \mathcal{E} \iff (\nu^j,\nu^i)\in \mathcal{E}$. An edge, denoted as $e^m:=(\nu^i,\nu^j)$, indicates that agent $i$ receives information from agent $j$, where $m$ represents the number of edge $(\nu^i,\nu^j)$. Let us enumerate the edge set as $\mathcal{E}=\{e^1,...,e^M\}$, where $e^m\in\mathcal{E}$ represents the $m$-th edge. For $m\in \{1,...,M\}$, let $\alpha^m\in \mathbb{R}$ be a positive scalar denoting the edge weight corresponding to the $m$-th edge.

The set of neighbors of agent $i$ is denoted by ${\cal N}^i$. The (oriented) incidence matrix $D \in \mathbb{R}^{N\times M}$ of the graph $\mathcal{G}$ is defined component-wise by
\begin{align*}
    D_{im} = \begin{cases}
    +1, & \text{if}\ \text{node}\ i\ \text{is}\ \text{the}\ \text{source}\ \text{node}\ \text{of}\ \text{edge}\ e^m,\\
    -1, & \text{if}\ \text{node}\ i\ \text{is}\ \text{the}\ \text{sink}\ \text{node}\ \text{of}\ \text{edge}\ e^m,\\
    0, & \text{otherwise},
    \end{cases}
\end{align*}
where for undirected graphs the orientation in the incidence matrix can be chosen arbitrarily.

The weighted Laplacian of a graph $\mathcal{G}$ can be defined as 
\begin{align*}
    L = DWD^T,
\end{align*}
where $W=\text{diag}(\alpha^1,...,\alpha^M)\in \mathbb{R}^{M\times M}$ is a diagonal matrix of edge weights. The Laplacian matrix is symmetric and positive semi-definite. In the context of game theory, we define  a local Laplacian for  agent / player $i$ as
\begin{align*}
    L^i = DW^iD^T,
\end{align*}
and let $W^i \in \mathbb{R}^{M\times M}$ be a diagonal matrix such that the $m$-th diagonal entry of $W^i$ is equal to $\alpha^m$ if $e^m\in \mathcal{E}^i$ and zero otherwise, where $\mathcal{E}^i = \{ e^i,...,e^i_{deg_i}\}\subset \mathcal{E}$ be the set of edges incident at node $\nu^i\in V$. In this paper we assume $\alpha^m = 1$, for all $\forall m\in M$.

\begin{assumption}
Graph $\mathcal{G}$ is connected, i.e.\  there exists an undirected path between every two vertices $\nu^i, \; \nu^j\in {\cal V}$,  $j\neq i$.
\end{assumption}
From now on, we assume that the graph used in this paper is undirected. 
% Let $X\in \mathbb{R}^{a\times b}$, $Y\in \mathbb{R}^{c\times d}$. The Kronecker product is defined as
% $$
% X\otimes Y = \begin{bmatrix}
% x_{11}Y & \cdots & x_{1b}Y\\
% \vdots & \ddots & \vdots\\
% x_{a1}Y & \cdots & x_{ab}Y
% \end{bmatrix}\in \mathbb{R}^{ac\times bd}.
% $$

\subsection{$\sigma$-Norms}
The $\sigma$-norm of a vector is a map $\mathbb{R}^n\rightarrow \mathbb{R}_{\geq 0}$ \textit{(not a norm)} defined as \cite{olfati}
\begin{align}\label{sigma}
    ||y||_\sigma = \frac{1}{\epsilon}\left[ \sqrt{1+\epsilon ||y||^2} - 1\right],
\end{align}
where $||\cdot||$ is an Euclidian norm in $\mathbb{R}^n$, $\epsilon >0$ is a small scalar value, and the gradient $\sigma_\epsilon(y) = \nabla ||y||_\sigma$ is
\begin{align*}
    \sigma_\epsilon(y) = \frac{y}{\sqrt{1+\epsilon ||y||^2}} = \frac{y}{1+\epsilon ||y||_\sigma}.
\end{align*}
The map $||y||_\sigma$ is differentiable everywhere. This property of the $\sigma$-norm will be used when dealing with the norm in the state-dependent weight matrix.

\subsection{Agent Dynamics}

In this article, we consider a homogeneous multi-agent system where each agent is modeled as a zero-order hold discretisation of a double integrator. Each agent is assumed to be moving in an $n$-dimensional plane. In the context of game theory, each agent acts as a player in the game. The single-agent discrete-time dynamics is
\begin{align}\label{dyn1}
    x^i_{k+1} = f x^i_k + g u^i_k, \hspace{3pc} \text{for}\ i = 1,...,N,
\end{align}
where the state vector for agent $i$ is $x^i_k = \begin{bmatrix} p^i_k, v^i_k\end{bmatrix}^T \in \mathbb{R}^{2n}$, and contains  position $p^i_k$ and velocity $v^i_k$ at time $k$, with
\begin{align*}
    f = \begin{bmatrix} 1 & \delta\\0 & 1\end{bmatrix}\otimes I_n \in \mathbb{R}^{2n\times 2n},\ \
    g = \begin{bmatrix} \frac{\delta^2}{2}\\ \delta\end{bmatrix}\otimes I_n \in \mathbb{R}^{2n\times n}.
\end{align*}
Here, $u^i_k$ is the (acceleration) control input of agent $i$, and $\delta$ is the sampling time.  To define the state vector for the multi-agent system, we select  $x_k = \begin{bmatrix} p^1_k,\cdots, p^N_k, 1, v^1_k, \cdots, v^N_k\end{bmatrix}^T \in \mathbb{R}^{2Nn+1}$. Having an entry with value $1$ between the positions and velocities allows the inclusion of a formation offset term, as explained below.  The multi-agent dynamics can then be represented as
\begin{align}\label{dyn2}
    x_{k+1} = F x_k + \sum_{i=1}^N G^i u^i_k, 
\end{align}
where
\begin{align*}
    F &= \begin{bmatrix} I_{Nn} & 0_{Nn\times 1} & \delta I_{Nn}\\
    0_{1\times Nn} & 1 & 0_{1\times Nn}\\
    0_{Nn\times Nn} & 0_{Nn\times 1 } & I_{Nn}
    \end{bmatrix} \in \mathbb{R}^{(2Nn+1)\times (2Nn+1)},
    \\
    G^i &= \begin{bmatrix} \frac{\delta^2}{2} \hat{g}^i\\ 0_{1\times n}\\ \delta \hat{g}^i \end{bmatrix} \in \mathbb{R}^{(2Nn+1)\times n},
    \end{align*} 
with $\hat{g}^i = \hat{c}^i\otimes I_n \in \mathbb{R}^{Nn\times n}$, where $\hat{c}^i$ is the $i$-th column of the identity matrix of size $N$. $I_{Nn} \in \mathbb{R}^{Nn\times Nn}$ is an identity matrix.  The scalar value of 1 in the matrix $F$ corresponds to a formation offset term, which will be explained in the next subsection.

\subsection{Formation with Collision Avoidance on the Nodes System}

The problem considered in this article is  formation control, i.e.\ all agents in a multi-agent system are supposed to move from arbitrary initial states to attain a formation (specified in terms of desired displacements $d^{ij}$ between agents $i$ and $j$), while minimizing a performance index over a finite time horizon $[0, T]$
\begin{align}\label{cost2}
    J^i(U^i) = \frac{1}{2} \Big( X_k^T \mathcal{Q}^i(x_k) X_k +  U_k^{i^T} \mathcal{R}^{ii} U^i_k \Big),
\end{align}
% \begin{align}\label{cost1}
%     J^i(U^1,..,U^N) = \frac{1}{2}\big( X_k^T \mathcal{Q}^i(x_k) X_k + \sum_{j=1}^N U_k^{j^T} \mathcal{R}^{ij} U^j_k\big),
% \end{align}
with the stacked state vector for the whole horizon $X_k = \begin{bmatrix} x_{k+1},x_{k+2},...,x_{k+T}\end{bmatrix}^T \in \mathbb{R}^{(2Nn+1)T}$ and the stacked control inputs vector 
$U^i_k = \begin{bmatrix} u^i_{k},u^i_{k+1},...,u^i_{k+T-1}\end{bmatrix}^T \in  \mathbb{R}^{NnT}$. The state weighting matrix for each agent $i$ is given by $\mathcal{Q}^i(x_k) = \text{blkdiag}(Q^i(x_k),...,Q^i(x_k),Q^i_T(x_T))\in  \mathbb{R}^{(2Nn+1)T\times (2Nn+1)T}$, where $Q^i(x_k) = (Q^i_\alpha + Q^i_\beta(x_k)) \in \mathbb{R}^{(2Nn+1)\times (2Nn+1)}$ is a positive semi definite matrix, with $Q^i_\alpha$ and $Q^i_\beta(x_k)$  represent the weighting matrices for formation and collision avoidance terms, respectively. The terminal weighting matrix $Q^i_T(x_T)$ has the same pattern as $Q^i(x_k)$ and can be defined by choosing arbitrary scalar weights of ${\beta}^i>0$.

The control weighting matrix is $\mathcal{R}^{ii} = \text{blkdiag}(R^{ii}) \in \mathbb{R}^{NnT\times NnT}$, where $R^{ii} \in \mathbb{R}^{Nn\times Nn}$ is a positive definite matrix. Here, we assume there is no cross coupling in the input, i.e., $\mathcal{R}^{ij} = 0$, where $j\neq i$. Next, the rest of this subsection is dedicated to discussing the formulation of the first term of the cost in \eqref{cost2}.
% Due to this reason, the individual cost functions in \eqref{cost1} is reduced to
The formation error of each agent $i$ with collision avoidance can be expressed as
\begin{align}\label{localerr}
   \Psi^i_k =  &\sum_{j\in \mathcal{N}^i} \Big\{ \big( || p^i_k - p^j_k - d^{ij} ||^2 + ||v^i_k - v^j_k||^2\big) \nonumber\\
    &\hspace{1pc}+{\beta}^i \Big(\frac{ || p^i_k - p^j_k - d^{ij} ||^2 + ||v^i_k - v^j_k||^2}{||p^i_k-p^j_k||^2 - r^{i^2}}\Big) \Big\},\end{align}   
where $r^i$ is the safety radius of agent $i$ that is assumed to be the same for all $i\in N$ homogeneous agent, i.e $r^i = r$, and $\beta^i>0$ is a tuning parameter for agent $i$. By the property of sum-of-squares, \eqref{localerr} can be transformed into a matrix form
\begin{small}
\begin{align*}
    &\sum_{j\in \mathcal{N}^i} \Big\{ || p^i_k - p^j_k ||^2 - 2(p^i_k-p^j_k)^Td^{ij}+ ||d^{ij}||^2 + ||v^i_k - v^j_k||^2 \\
    &\hspace{1.5pc}+  \frac{\beta^i || p^i_k - p^j_k ||^2}{||p^i_k-p^j_k||^2 - r^2} - \frac{2\beta^i (p^i_k-p^j_k)^Td^{ij}}{||p^i_k-p^j_k||^2 - r^2}\\
    &\hspace{3pc}+  \frac{\beta^i ||d^{ij}||^2}{||p^i_k-p^j_k||^2 - r^2} +\frac{\beta^i || v^i_k - v^j_k ||^2}{||p^i_k-p^j_k||^2 - r^2}  \Big\}=
\end{align*}\end{small}
\vspace{-1pc}
\begin{align*}
    &{p_k^T\mathcal{L}^i_\alpha p_k - 2p_k^T\mathcal{D}\mathcal{W}^i_\alpha d + d^T\mathcal{W}^i_\alpha d + v_k^T \mathcal{L}^i_\alpha v_k + p_k^T\mathcal{L}^i_\beta(x_k) p_k} \\
    &\hspace{1pc}  - 2p_k^T\mathcal{D}\mathcal{W}^i_\beta(x_k) d + d^T\mathcal{W}^i_\beta(x_k) d + v_k^T \mathcal{L}^i_\beta(x_k) v_k\\
    &\hspace{2pc}= x_k^T\big(Q^i_\alpha + Q^i_\beta(x_k)\big)x_k= x^T_k Q^i(x_k) x_k
\end{align*}
where
\begin{align*}
    Q^i_\alpha = \delta \begin{bmatrix}
    \mathcal{L}^i_\alpha & -\mathcal{D}\mathcal{W}^i_\alpha d & 0\\
    -(\mathcal{D}\mathcal{W}^i_\alpha d)^T & d^T\mathcal{W}^i_\alpha d & 0\\
    0 & 0 & \mathcal{L}^i_\alpha 
    \end{bmatrix} \end{align*}
has size $ \mathbb{R}^{(2Nn+1)\times (2Nn+1)}$, with a diagonal matrix with the edge weight $\mathcal{W}^i_\alpha = W^i \otimes I_n\in \mathbb{R}^{Mn\times Mn}$. A lifted local Laplacian matrix is defined as $\mathcal{L}^i_\alpha = \mathcal{D}\mathcal{W}^i_\alpha \mathcal{D}^T \in \mathbb{R}^{Nn\times Nn}$ with  $\mathcal{D} = D\otimes I_n \in \mathbb{R}^{Nn\times Mn}$ being the incidence matrix lifted to dimension $n$ of the space in which agents are moving, and $d=\text{col}(d^{ij}) \in \mathbb{R}^{Mn}$  the column vector of desired displacements vector $d^{ij}\in \mathbb{R}^{n}$. The state-dependent weighting matrix is then
\begin{align*}
    Q^i_\beta(x_k)   = \delta \begin{bmatrix}
    \mathcal{L}^i_\beta(x_k) & -\mathcal{D}\mathcal{W}^i_\beta(x_k) d & 0\\
    -(\mathcal{D}\mathcal{W}^i_\beta(x_k) d)^T & d^T\mathcal{W}^i_\beta(x_k) d & 0\\
    0 & 0 & \mathcal{L}^i_\beta(x_k)   \end{bmatrix} \end{align*}
of size $\mathbb{R}^{(2Nn+1)\times (2Nn+1)}$, with the state-dependent Laplacian matrix  defined as $\mathcal{L}^i_\beta(x_k) = \mathcal{D}\mathcal{W}^i_\beta(x_k) \mathcal{D}^T \in \mathbb{R}^{Nn\times Nn}$. The state-dependent edge weight matrix is $\mathcal{W}^i_\beta(x_k) = W^i_\beta(x_k)\otimes I_n\in \mathbb{R}^{Mn\times Mn}$,
% \begin{align*}
%     \mathcal{W}^i_\beta(x_k) = \text{diag}\Big(0,\cdots, \frac{\beta^i}{||p^i_k-p^j_k||^2 - r^2},\cdots,0\Big)\otimes I_n
% \end{align*}
% of size $\mathbb{R}^{Mn\times Mn}$
where now the $m$-th diagonal entry of $W^i_\beta(x_k)\in \mathbb{R}^{M\times M}$ is equal to  $\dfrac{\beta^i}{||p^i_k-p^j_k||^2-r^2}$ if $e^i\in \mathcal{E}^i$ and zero otherwise. The Laplacian matrix $\mathcal{L}^i_\beta(x_k)$ depends on the state since the diagonal edge matrix $\mathcal{W}^i_\beta(x_k)$ contains collision terms between agents $i$ and $j$. 

The formulation of the state vector $x_k \in \mathbb{R}^{2Nn+1}$ has been confirmed, and as a result, the state matrix $F$ matches the dimensions of the state weighting matrix $Q^i(x_k) \in \mathbb{R}^{(2Nn+1)\times (2Nn+1)}$. 

\begin{assumption}
The initial positions of the agents satisfy $\| p^i_0 - p^j_0\| > r^i+r^j$, for all $i,j\in N$, $j\neq i$.
\end{assumption}

By adopting the same reasoning as outlined in \cite{thulasi}, assumption 1 ensures that the term $x_0^T Q^i(x_0) x_0$ in \eqref{cost2}, for all $i \in N$, remains bounded. It follows that the agents operate without entering the avoidance region.

% As for the the second term in \eqref{localerr}, here we assume that the collision avoidance mechanism is capable of ensuring that the safety radius is not violated.

% implies that not only must the formation be maintained, but the distance between neighboring agents must also exceed the safety radius to prevent collisions. Thus, an assumption can be made as follows.
% \begin{assumption}
% The term $||p^i_k-p^j_k||^2 > r^2$ holds.
% \end{assumption}

\subsection{Nash Equilibrium and Coupled State-Dependent Riccati Equation (CSDRDE)}

The formulation of the formation control problem with dynamics  \eqref{dyn2} and cost functions  \eqref{cost2} as a game reflects the non-cooperative behavior, where each player is searching for a Nash equilibrium corresponding to its own local cost function. 
\begin{definition}
A collection of strategies $U^{i\star}$ constitutes a Nash equilibrium if and only if the inequalities
    $$
    J^i(U^{1\star},...,U^{N\star}) \leq J^i(U^{1\star},...,U^{i-1\star},U^i,U^{i+1\star},...,U^{N\star})
    $$ 
hold for $i=1,...,N$.    
\end{definition}
We now  formulate the first problem (for the multi-agent system running on the nodes) as follows.

\begin{problem}\label{prob1}
Find  local control sequences that achieve a Nash equilibrium corresponding to the local cost functions \eqref{cost2} over the control input sequences $u^i$ subject to \eqref{dyn2}.
\end{problem}

\begin{theorem}
An open-loop Nash equilibrium for the game defined by Problem \ref{prob1} is achieved by the control sequences
\begin{align}\label{unash}
    u^{i\star}_k(x_k) =  K^i_k(x_k) x_k,
\end{align}
where
\begin{align}\label{Nashgain}
    K^i_k(x_k) &= -R^{{ii}^{-1}} G^{i^T} P^i_{k+1}(x_{k+1})\Lambda_k^{-1}F,
\end{align}
and $P^i_{k+1}(x_{k+1})$ is the solution to the coupled state-dependent Riccati difference equation
\begin{align}\label{sdre}
    P^i_k(x_k) &= F^T P^i_{k+1}(x_{k+1})\Lambda_{k+1}^{-1}F+ Q^i(x_k) \nonumber\\
    &\hspace{-0.5pc}+(I_N \otimes x_k^T)\begin{bmatrix}
		x_k^T\frac{\partial Q^i(x_k)}{\partial x^1_k},&...,&x_k^T\frac{\partial Q^i(x_k)}{\partial x^N_k}
		\end{bmatrix}^T,
\end{align}
which can be solved backward with $P^i_{T}(x_T) = Q^i_T(x_T)$. The corresponding closed-loop state trajectory is 
\begin{align}\label{closedloopNash}
    x_{k+1}^\star = \Lambda_k^{-1}Fx_k^\star,
\end{align}
where 
\begin{align}\label{lamb}
    \Lambda_k = \Big( I + \sum_{j=1}^N G^jR^{{jj}^{-1}}G^{j^T}P^j_{k+1}(x_{k+1})\Big).
\end{align}
\end{theorem}

\begin{proof}
See appendix.
\end{proof}
By looking at \eqref{sdre}, solving for $P_k^i(x_k)$ requires information about $P_{k+1}^j(x_{k+1})$ for all $j\in N$, and thus the Riccati equation cannot be solved in a distributed way. Furthermore, it should be noted that once the Riccati equation is solved, the state feedback gains obtained in \eqref{Nashgain} are fully populated and require knowledge of all states across the network.

\section{DISTRIBUTED FRAMEWORK}
% As mentioned in the previous section, the original problem cannot be solved in a distributed manner due to the coupling terms in the cost functions when solving on the nodes. Therefore, 

Building upon the method from \cite{Prima}, this section outlines a distributed strategy to address the issue. The strategy involves an associated fictitious multi-agent system that evolves on the edges of the communication graph, departing from the conventional node-based approach.

\subsection{The Edge System}
Inspired by \cite{Jond2}, we associate a fictitious agent with each edge ($\nu^i,\nu^j$) of the communication graph with dynamics
\begin{align}\label{newstate}
\begin{bmatrix}q^m_k\\w^m_k \end{bmatrix} = 
\begin{bmatrix} p^i_k - p^j_k - d^{ij}\\
 v^i_k - v^j_k \end{bmatrix} \hspace{1pc}\text{and}\hspace{0.5pc} a^m_k = u^i_k - u^j_k, \hspace{1pc} 
\end{align}
for $m=1,...,M$. The state vector for  edge agent $m$ is  $z^m_k = \begin{bmatrix}
q^m_k, w^m_k\end{bmatrix}\in \mathbb{R}^{2n}$.
Then, the  relative dynamics for  edge agent $m$ is 
\begin{align*}
    z^m_{k+1} = fz^m_k + g a^m_k,\hspace{1pc}\text{for}\ m=1,.,,,M.
\end{align*}
The state vector for the whole edge system can be arranged as $\Tilde{z}_k = [z_k^1,...,z_k^M]^T\in \mathbb{R}^{2Mn}$. We rearrange the states by a permutation
\begin{align*}
     z_k = \Pi \Tilde{z}_k,
\end{align*}
with permutation matrix 
 \begin{align*}
     \Pi = \begin{bmatrix}
      I_M\otimes \begin{bmatrix}
       1&0\end{bmatrix}\\
       I_M\otimes \begin{bmatrix}
       0&1\end{bmatrix}
     \end{bmatrix}\otimes I_n \in \mathbb{R}^{2Mn\times 2Mn}.
 \end{align*}
Therefore, the whole edge dynamics can be written as
\begin{align}\label{dyn3}
    z_{k+1} = \bar{F} z_k + \sum_{m=1}^M \bar{G}^m a^m_k,
\end{align}
where $z_k = [q^1_k,...,q^M_k,w^1_k,...,w^M_k]^T\in \mathbb{R}^{2Mn}$ and  
\begin{align*}
    \bar{F} &= \left(\begin{bmatrix} 1 & \delta \\ 0 & 1 \end{bmatrix}\otimes I_M \otimes I_n\right)\in \mathbb{R}^{2Mn\times 2Mn}, \hspace{1pc} \\
    \bar{G}^m &= \begin{bmatrix} \frac{\delta^2}{2} \bar{g}^m \\ \delta \bar{g}^m \end{bmatrix} \in \mathbb{R}^{2Mn\times n},
    \end{align*} 
with $\bar{g}^m = \bar{c}^m\otimes I_n \in \mathbb{R}^{Mn\times n}$, where $\bar{c}^m$ is the $m$-th column of identity matrix of size $M$.

\subsection{Formation with Collision Avoidance on the Edge System}

Since we  relocated the coupling terms that were initially in the cost function to the system dynamics,  the local error for an edge agent $m$ at time instance $k$ to be minimized is
\begin{align*}
   \bar{\Psi}^m_k &=    \alpha^m  \big( || q_k^m ||^2 + ||w_k^m||^2\big) +
    {\beta}^m \Big(\frac{ || q^m_k ||^2 + ||w^m_k||^2 }{||q^m_k + d^{ij}||^2 - r^2}\Big)\\
    &= z_k^T \big( \bar{Q}^m_\alpha + \bar{Q}^{m}_\beta(z_k^m)\big) z_k\nonumber\\
    &= z_k^T \bar{Q}^{m}(z_k^m) z_k\nonumber
\end{align*}
where
\begin{align*}
    \bar{Q}^{m}_\alpha &= \delta ( I_2 \otimes \bar{W}^{m}_\alpha \otimes I_n) \in \mathbb{R}^{2Mn\times 2Mn},\\ 
    % \ \text{with}\\ \bar{W}^{m}_\alpha &= \text{diag}(0,\cdots,\alpha^m,\cdots,0)\in \mathbb{R}^{M\times M}, 
% \end{align*}
% and
% \begin{align*}
    \bar{Q}^{m}_\beta(z_k^m) &= \delta ( I_2 \otimes \bar{W}^{m}_\beta(z_k^m) \otimes I_n) \in \mathbb{R}^{2Mn\times 2Mn},
    % \hspace{0.5pc}\ \text{with}\\ \bar{W}^{m}_\beta(z_k^m) &= \text{diag}\Big(0,...,\frac{{\beta}^m}{||q^m_k + d^{ij}||^2 - r^2},...,0\Big)\in \mathbb{R}^{M\times M} .
\end{align*}
with $\bar{W}^{m}_\alpha\in \mathbb{R}^{M\times M}$ is a diagonal matrix such that the $m$-th diagonal entry of $\bar{W}^{m}_\alpha$ is equal to $\alpha^m$ and zero otherwise and let $\bar{W}^{m}_\beta(z_k^m)\in \mathbb{R}^{M\times M}$ be a diagonal matrix such that the $m$-th diagonal entry of $\bar{W}^{m}_\beta(z_k^m)$ is equal to $\dfrac{{\beta}^m}{||q^m_k + d^{ij}||^2 - r^2}$ and zero otherwise.

Note that in contrast to $Q^i(x_k)$ from the first problem, $\bar{Q}^m(z^m_k)$ here is a block diagonal matrix where edge dynamics are decoupled. Therefore, we can arrange the decoupled cost function for the $m$-th edge as 
    \begin{align}\label{cost3}
    \bar{J}^{m}(A^m) = \frac{1}{2}\Big( Z_k^T\mathcal{\bar{Q}}^m(z^m_k)Z_k +  A^{m^T}_k\mathcal{\bar{R}}^{mm}A^m_k\Big),
\end{align}
where the stacked edge state vector now is arranged as $Z_k = [z_{k+1},z_{k+2},...,z_{k+T}]^T\in \mathbb{R}^{2MnT}$ and the stacked relative control inputs vector is $A^m_k = [a_k^m,a_{k+1}^m,...,a_{k+T-1}^m]^T\in \mathbb{R}^{MnT}$. 

The state  weighting matrix for the new cost evolving on edges is defined as $\mathcal{\bar{Q}}^m(z^m_k) = \text{blkdiag}\big(\bar{Q}^m(z_k^m),...,\bar{Q}^m(z_k^m),\bar{Q}_T^m(z_T^m)\big)\in \mathbb{R}^{2MnT\times 2MnT}$, where the terminal cost $\bar{Q}_T^m(z_T^m)\in \mathbb{R}^{2MnT\times 2MnT}$ has the same pattern as $\bar{Q}^m(z_k^m)$ with arbitrary choices of scalar weights instead of $\alpha^m, \beta^m > 0$. The control weight is $\mathcal{\bar{R}}^{mm} = \text{blkdiag}(\bar{R}^{mm}) \in \mathbb{R}^{MnT\times MnT}$ with a positive definite matrix $\bar{R}^{mm}\in \mathbb{R}^{Mn\times Mn}$. Finally, we  formulate the new problem for the  edge dynamics \eqref{dyn2} as follows.
\begin{problem}\label{prob2}
Minimize the local cost function \eqref{cost3} over the relative acceleration control input sequences $a^i$ subject to dynamics \eqref{dyn3}.
\end{problem}

\begin{theorem}
The optimal solution to Problem \ref{prob2} is 
\begin{align}\label{optimalcont}
    a^{m\star}_k(z_k^m) = \bar{K}^m_k(z_k^m) z_k,\hspace{1pc}\text{for} \ m = 1,...,M,
\end{align}
where
\begin{align}\label{gainK}
    \bar{K}^m_k(z_k^m) &= -(\bar{R}^{mm} + \bar{G}^{m^T}\bar{P}^m_{k+1}(z_{k+1}^m)\bar{G}^m)^{-1}\times\nonumber\\
    &\hspace{7pc}\bar{G}^{m^T}\bar{P}^m_{k+1}(z_{k+1}^m)\bar{F},
\end{align}
and $\bar{P}^m_{k+1}(z_{k+1}^m)$ is the solution to the decoupled state-dependent Riccati difference equation
\begin{align}\label{dsdrde}
    &\bar{P}^m_k(z_k^m) = \bar{F}^T\bar{P}^m_{k+1}(z_{k+1}^m)\bar{F}+ \bar{F}^T\bar{P}^m_{k+1}(z_{k+1}^m)\bar{G}^m \bar{K}^m_k(z_k^m) \nonumber\\
			&\hspace{-0.4pc}+ \bar{Q}^m(z_k^m) + (I_M\otimes z^{m^T}_k) \begin{bmatrix}z_k^{m^T}
			\frac{\partial \bar{Q}^m(z_k^m)}{\partial z^1_k}\ ...\ z_k^{m^T}\frac{\partial \bar{Q}^m(z_k^m)}{\partial z^M_k}
			\end{bmatrix}^T 
\end{align}
with $\bar{P}^m_T(z_T^m) = \bar{Q}_T^m(z_T^m)$.
\end{theorem}

\begin{proof}
Can be shown similarly to Section 3 in \cite{Chang}.
\end{proof}
Note that because both $\bar{Q}^m(z_k^m)$ and its derivative in \eqref{dsdrde} involve the norm of a variable, the $\sigma$-norm defined in \eqref{sigma} is employed to ensure differentiability throughout. The feedback gains $\bar{K}_k^m(z_k^m)$ in \eqref{gainK} are now decoupled from each other. This decoupling principally permits a distributed implementation, in contrast to $K_k^i(x_k)$ in \eqref{Nashgain}.

However, solving the decoupled SDRDE in \eqref{dsdrde} is challenging due to its state-dependency. To address this challenge, we embrace the receding horizon technique for solving SDRDE presented in \cite{Dutka}. The approach involving the decoupled SDRDE entails the following steps:
\begin{enumerate}
    \item[1.] Utilize the state-feedback gains $\bar{K}^m_k$ from \eqref{gainK} computed in the previous iteration. Let $z^{p}_k$ be the prediction of the dynamics, commencing from the current state $z_k$.
    \item[2.] Work backwards in time to compute the Riccati solution, yielding $\bar{P}^m_{k+T},...,\bar{P}^m_{k+1}$ along the predicted state trajectory.
    \item[3.] Employ this information to update the state feedback gains $\bar{K}^m_k,...,\bar{K}^m_{k+T-1}$. Implement the first gain $\bar{K}^m_k$ for control purposes.
    \item[4.] At the subsequent sampling instant, repeat this process, and make use of the remaining gains $\bar{K}^m_{k+1},...,\bar{K}^m_{k+T-1}$.
    \item[5.] Determine the terminal gain $\bar{K}^m_{k+T}$ required for the next iteration by solving the decoupled SDRDE along the predicted states. This approach facilitates a receding horizon strategy.
\end{enumerate}
The detailed steps to evaluate the decoupled SDRDE approach are provided in Algorithm 1.
\begin{algorithm}[!h]\label{algo1}
    \caption{DSDRDE \cite{Dutka}}
  \begin{algorithmic}[1]
        \INPUT $\bar{Q}^m_\alpha$, $\bar{Q}^m_\beta(z^m_k)$ and its derivative $\frac{\partial \bar{Q}^m_\beta(z^m_k)}{\partial z^m}$, horizon $T$, number of edges $M$, tolerance $\varepsilon$, and $t_{max}$
        \OUTPUT The relative control inputs $a_k^{m\star}$
        \STATE At  time $k$ do the following
        \STATE $t = 1$ 
        \STATE Initialize $z^{p}_{t=1} = z_k$ 
        \STATE Initialize $\bar{K}^m = \text{zeros}(M,2M,T)$ 
        \WHILE {$t\leq t_{max}$ and $||\bar{K}_{old}^m-\bar{K}^m||>\varepsilon$}
            \STATE $\bar{K}_{old}^m = \bar{K}^m$
            \FOR {$j=1:T$}
            \STATE $z^{p}_{j+1} = \Big( \bar{F} + \sum_{m=1}^M\bar{G}^m\bar{K}^m\Big) z^{p}_{j}$
            \ENDFOR
            \STATE Set $ \bar{P}_T^m = \bar{Q}^m_\beta(z^{p}_T)$
            \FOR {$j = T:-1:2$}
            \STATE Calculate $\bar{P}^m_{j-1}$  from \eqref{dsdrde}
            \STATE Calculate $\bar{K}^m(:,:,j-1)$ from \eqref{gainK}
            \ENDFOR
            \STATE $t \leftarrow t+1$
        \ENDWHILE
        \STATE Set $a_k^{m\star} = \bar{K}^m(:,:,j)z_k$ 
    \end{algorithmic}
\end{algorithm}

% After obtaining the optimal solution on the edge system, i.e., $a_k^{m\star}$, the question of how to get back the optimal solution on the nodes arises.

\subsection{Distributed Implementation}
In this subsection, we show how to obtain the optimal control inputs of the physical vertex agents from the relative control inputs $a^{m\star}_k$ of the fictitious (edge) agents in a distributed fashion.  We will use the symbol $\hat{u}^{i\star}_k$ to denote the   physical control inputs corresponding to the fictitious relative control inputs $a^{m\star}_k$. Recall that from \eqref{newstate}, we can express the relation between $a^{\star}_k$ and $\hat{u}^{\star}_k$ as
\begin{align*}
    \Phi \hat{u}^\star_k =  a^\star_k,
\end{align*}
where $\Phi = D^T\otimes I_n$ and $a^{\star}_k = [a^{{1\star}^T}_k,...,a^{{M\star}^T}_k]^T$.

We  consider  minimizing the residual $f({u}) = ||\Phi \hat{u}^\star_k -  a^\star_k||^2$. Since the undirected graph $\mathcal{G}$ is assumed to be connected,  there exists a unique solution to  minimizing the residual $f({u})$, given by 
\begin{align}
    \label{udagger}
    \hat{u}^\star_k = \Phi^\dagger a^\star_k,
\end{align}
where $\Phi^\dagger$ is the pseudo-inverse of $\Phi$. Since $\Phi^\dagger$ is a fully populated matrix,  this  will lead to a centralized solution. To compute \eqref{udagger} in a distributed way, a distributed steepest descent algorithm  is employed, which updates the local control input at iteration step $l$ according to 
\begin{align}\label{gs}
    \hat{u}^\star_{l+1} = (I - 2\gamma \Phi^T\Phi )\hat{u}^\star_l + 2\gamma \Phi^T a^\star_k,
\end{align}
with $\gamma$ as a learning rate that satisfies
$$ 
2\gamma \leq \frac{2}{\lambda_{\rm max}(\Phi^T\Phi)}.
$$ It was demonstrated in \cite{fbullo} that this algorithm converges to a solution $\hat{u}^\star_k$ in \eqref{udagger} which is unique. The key fact is that the two matrices on the right-hand side of \eqref{gs} are sparse and allow a distributed computation of the updates $\hat{u}^\star_{l+1}$. The detailed steps to evaluate this approach are provided in Algorithm 2.

\begin{algorithm}[!h]\label{algo2}
    \caption{Iterative Distributed Steepest Descent on Updating Local Control Inputs \cite{Prima}}
  \begin{algorithmic}[1]
        \INPUT Tuning parameter $\gamma$, tolerance $\varepsilon$, and $l_{max}$
        \STATE {At  time $k$ do the following}
        \STATE{$l=1$}
        \STATE {Initialize $\hat{u}^i_{l=1}$, for $i=1,...,N$ with \textit{"warm start"}}
        \WHILE{$l\leq l_{max}$ and $||\Phi \hat{u}_l-a^\star_k||>\varepsilon$} 
        \STATE \textbf{Receive} $\hat{u}^j_{l}$ \textbf{from agent(s)} $j \in \mathcal{N}^i$
        \STATE{Calculate $\hat{u}^i_{l+1} = \hat{u}^i_l - 2\gamma \sum_{j\in \mathcal{N}^i}\big( \hat{u}^i_l - \hat{u}^j_l - a^\star_{(ij)}\big)$}
        \STATE \textbf{Broadcast} $\hat{u}^i_{l+1}$ \textbf{to agent(s)} $j \in \mathcal{N}^i$
        \STATE{$l \leftarrow l+1$}
        \ENDWHILE
    \end{algorithmic}
\end{algorithm}

The sixth step in Algorithm 2 can be interpreted as follows: during each iteration, the estimate at node $i$ is updated based on the errors in the relative control inputs, and each error in the edge (i.e., the difference between the estimated and measured edge difference) contributes to a  correction in the node value.

\section{ILLUSTRATIVE EXAMPLE}

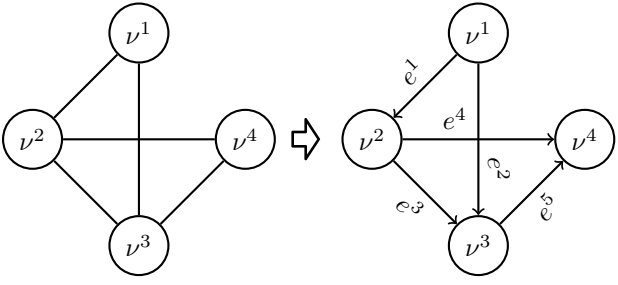
\begin{figure}[!tbhp]
\begin{minipage}{9.5pc}
    \centering
    \begin{tikzpicture}[node distance={20mm}, thick, main/.style = {draw, circle}] 
    \node[main] (1) {$\nu^2$}; 
    \node[main] (2) [above right of=1] {$\nu^1$}; 
    \node[main] (3) [below right of=1] {$\nu^3$}; 
    \node[main] (4) [above right of=3] {$\nu^4$}; 
    \draw[-] (2) -- node[midway, above left, sloped, pos=0.25] {} (1); 
    \draw[-] (1) -- node[midway, below right, sloped, pos=0.25] {} (3); 
    \draw[-] (3) -- node[midway, below left, sloped, pos=0.75] {} (4); 
    \draw[-] (2) -- node[midway, above right, sloped, pos=0.55] {} (3); 
    \draw[-] (1) -- node[midway, above right, sloped, pos=0.2] {} (4); 
    \end{tikzpicture}\end{minipage}$\vpointer$\hspace{0.3pc}
   \begin{minipage}{9pc}
    \begin{tikzpicture}[node distance={20mm}, thick, main/.style = {draw, circle}] 
    \node[main] (1) {$\nu^2$}; 
    \node[main] (2) [above right of=1] {$\nu^1$}; 
    \node[main] (3) [below right of=1] {$\nu^3$}; 
    \node[main] (4) [above right of=3] {$\nu^4$}; 
    \draw[->] (2) -- node[midway, above left, sloped, pos=0.25] {$e^1$} (1); 
    \draw[->] (1) -- node[midway, below right, sloped, pos=0.25] {$e^3$} (3); 
    \draw[->] (3) -- node[midway, below left, sloped, pos=0.75] {$e^5$} (4); 
    \draw[->] (2) -- node[midway, above right, sloped, pos=0.55] {$e^2$} (3); 
    \draw[->] (1) -- node[midway, above right, sloped, pos=0.2] {$e^4$} (4); 
    \end{tikzpicture} \end{minipage}
    \caption{ Arbitrary orientation of the $M=5$ edges of an undirected graph with $N=4$ nodes.}
    \label{graph1}
\end{figure}

This section illustrates the proposed approach with a formation control problem where double integrator agents are moving in $n = 2$ dimensional  space.  We consider $N=4$ agents and $M=5$ edges with an undirected communication graph,  as displayed in  Figure \ref{graph1}.

The incidence matrix  is 
\begin{align*}
    D = \begin{bmatrix}
     1 & 1 & 0 & 0 & 0\\
     -1 & 0 & 1 & 1 & 0\\
     0 & -1 & -1 & 0 & 1\\
     0 & 0 & 0 & -1 & -1
    \end{bmatrix}\in \mathbb{R}^{4\times 5}.
\end{align*}
We assume that all agents  have zero initial velocities, except agent 1 that has $v^1_0 = \begin{bmatrix}0.5, 1\end{bmatrix}^T$. The agents have initial positions 
\begin{align*}
    p^1_0 &= \begin{bmatrix}3.5\\ 1\end{bmatrix}, p^2_0 = \begin{bmatrix}12\\ 1\end{bmatrix}, p^3_0 = \begin{bmatrix}0\\ 5\end{bmatrix},  p^4_0 = \begin{bmatrix}15\\ 3.5\end{bmatrix},
\end{align*}
with the desired displacements vectors and safety radius
\begin{align*}
    d^{12} &= \begin{bmatrix}1.5\\ 1\end{bmatrix}, \hspace{1pc}d^{13} = \begin{bmatrix}0\\ 2\end{bmatrix}, \hspace{1pc}d^{23} = \begin{bmatrix}-1.5\\ 1\end{bmatrix}, \\
    d^{24} &= \begin{bmatrix}-3\\ 0\end{bmatrix}, \hspace{1pc} d^{34} = \begin{bmatrix}-1.5\\ -1\end{bmatrix}, \hspace{0.4pc} r = 0.5.
\end{align*}
\subsection{Simulation Results}
We first show the evolution of  agents' positions if collision avoidance is ignored, by taking $\bar{Q}^m_\beta(z_k^m) = 0$, for all $m\in M$.

\begin{figure}[thpb]
    \centering
    \includegraphics[width=0.38\textwidth]{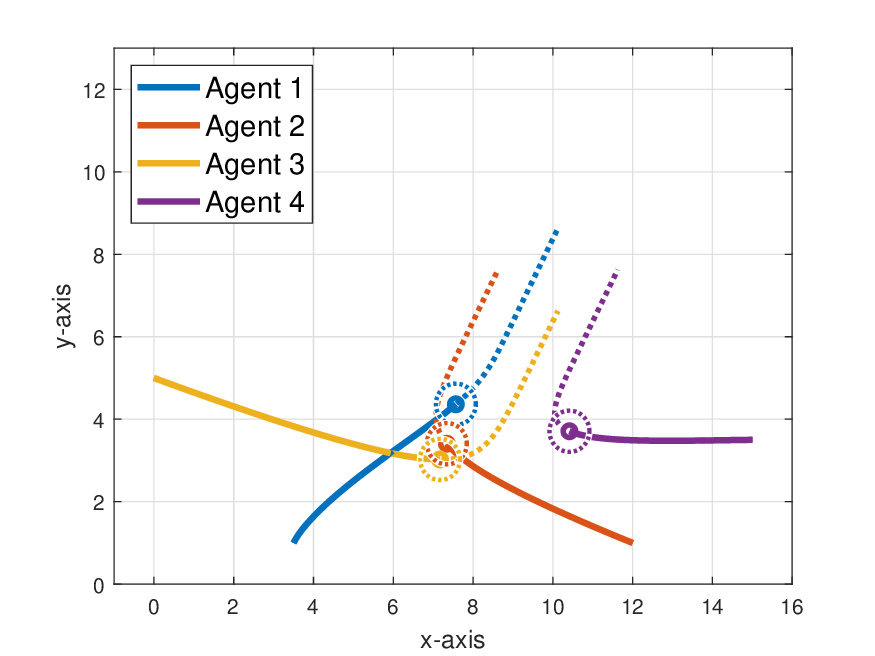}
    \caption{Progression of four agents' position on $x,y$-axes in 4 seconds, without collision avoidance.}
    \label{form0}
\end{figure}

Fig. \ref{form0} depicts the  four agents moving in the $x-y$-plane. A sampling time of 100$ms$ is used, and the dashed lines display the final trajectories over a period of 20$s$. The solid lines represent the intermediate progression run over 4$s$. The figure illustrates that the collision between agents within 4$s$.

For simulating formation control with collision avoidance, we construct $\bar{Q}^{m}_\beta(z_k^m)$ individually for each edge, wherein $\beta^m$ are set to 1 for all $m\in M$. Initially, we execute the steps associated with the DSDRDE using a horizon of $T = 10$, yielding the relative control inputs $a_k^{m\star}$. Once these relative inputs are acquired, we proceed with running the distributed steepest descent method to compute the physical control inputs $\hat{u}^i_k$ and to simulate the actual dynamics.

% We then use a horizon of $T = 10$, a threshold of $\varepsilon = 10^{-3}$, and $t_{max} = 12$ in Algorithm 1 to obtain the relative control inputs $a_k^{m\star}$. Finally, we execute Algorithm 2 to obtain the physical control inputs $\hat{u}^i_k$ and simulate the real dynamics by taking $\gamma = 0.1$ while setting $l_{max}= 10$.

\begin{figure}[!h]
\begin{minipage}{10pc}
    \centering
    \includegraphics[width=1.1\textwidth]{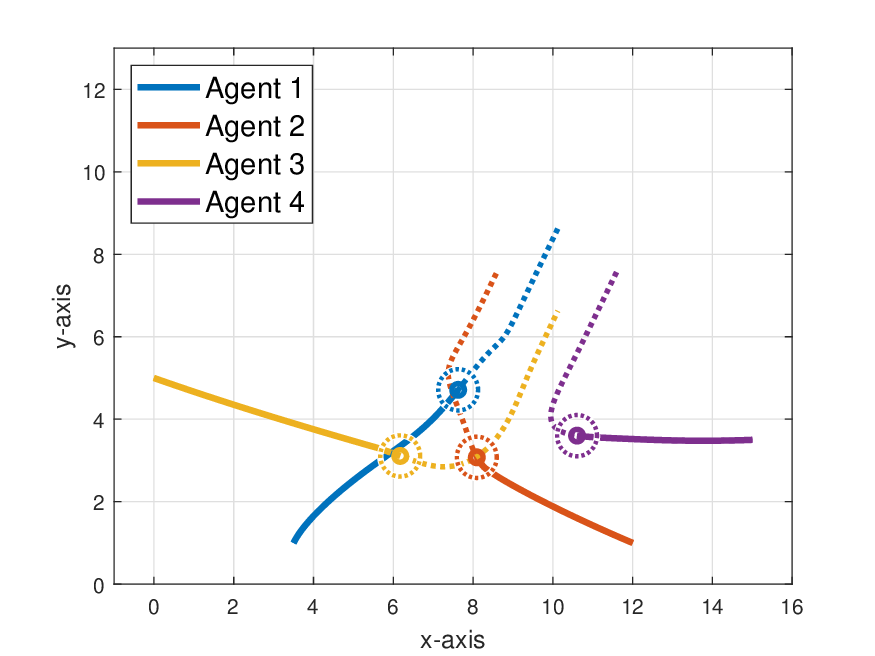}
    \caption{Progression of four agents' position on $x,y$-axes with  collision avoidance  in 4s.}
    \label{form4}
\end{minipage}
\begin{minipage}{10pc}
    \centering
    \includegraphics[width=1.1\textwidth]{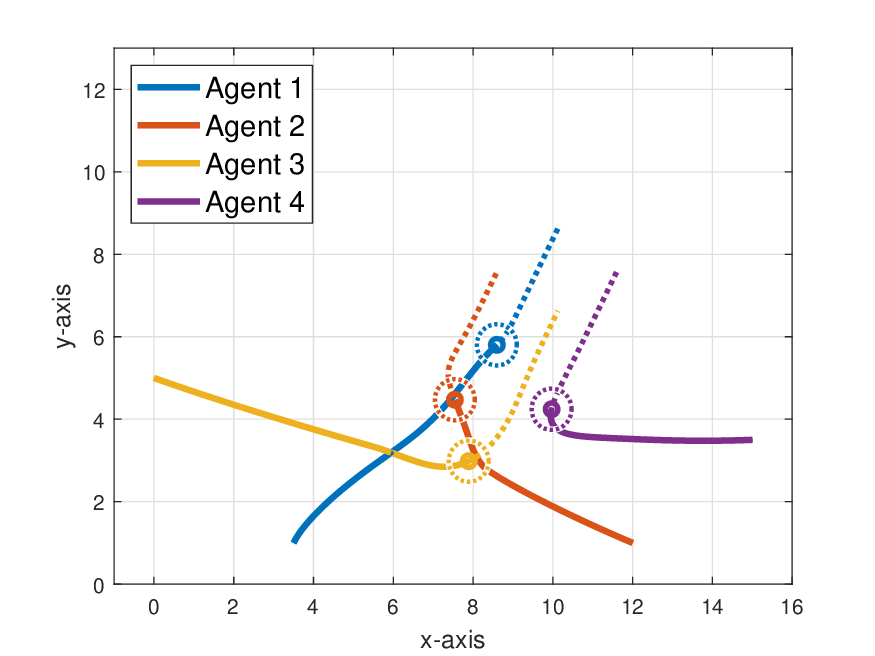}
    \caption{Progression of four agents' position on $x,y$-axes with  collision avoidance  in 7s.}
    \label{form7}
\end{minipage}      
\end{figure}

% \begin{figure}[!h]
%       \centering
%       \framebox{\parbox{3in}{
%   \includegraphics[width=0.44\textwidth]{figures/4sec.eps}}}
%         \caption{Progression of four agents' position on $x,y$-axes with  collision avoidance  in 4 seconds.}
%         \label{form4}
% \end{figure}

Formation with collision avoidance is visually presented in Figs. \ref{form4} and \ref{form7}. It is run for 4 and 7$s$, respectively. As depicted in Fig. \ref{form4}, agents one, two, and three successfully avoid collisions, in contrast to the scenario shown in Fig. \ref{form0}. When we extend the simulation time to 7$s$, a noteworthy observation emerges: agent three closely follows agent two, who, in turn, tracks agent one, resulting in their alignment within the desired formation. 

% \begin{figure}[thpb]
%       \centering
%       \framebox{\parbox{3in}{
%   \includegraphics[width=0.43\textwidth]{figures/7sec.eps}}}
%         \caption{Progression of four agents' position on $x,y$-axes with considering collision avoidance  in 7 seconds.}
%         \label{form7}
% \end{figure}

\begin{figure}[!h]
\begin{minipage}{10pc}
    \centering
    \includegraphics[width=1.1\textwidth]{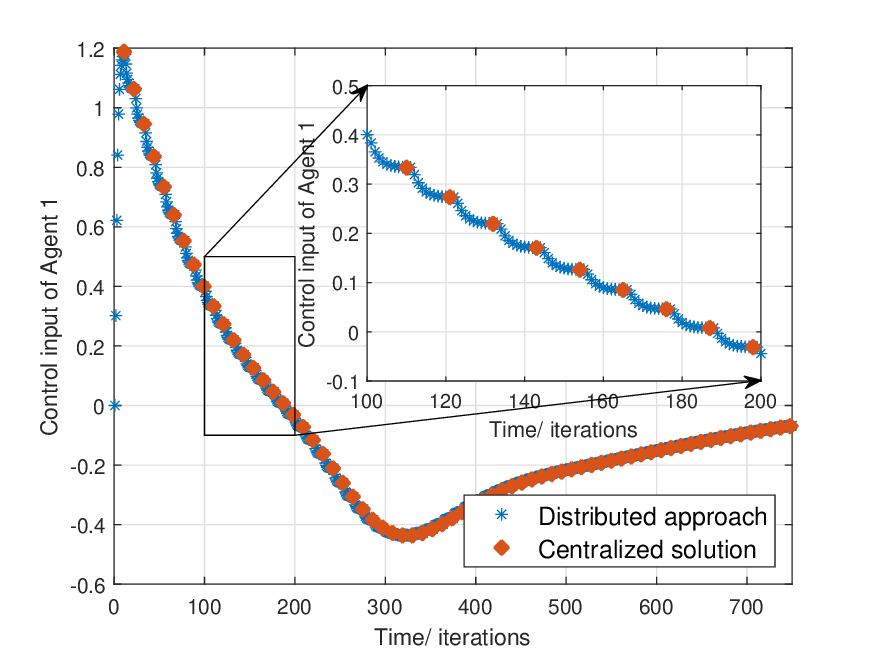}
    \caption{Control input of agent 1; centralised and distributed solution with 10 iterations/ sampling interval.}
        \label{continput1}
\end{minipage}
\begin{minipage}{10pc}
    \centering
    \includegraphics[width=1.1\textwidth]{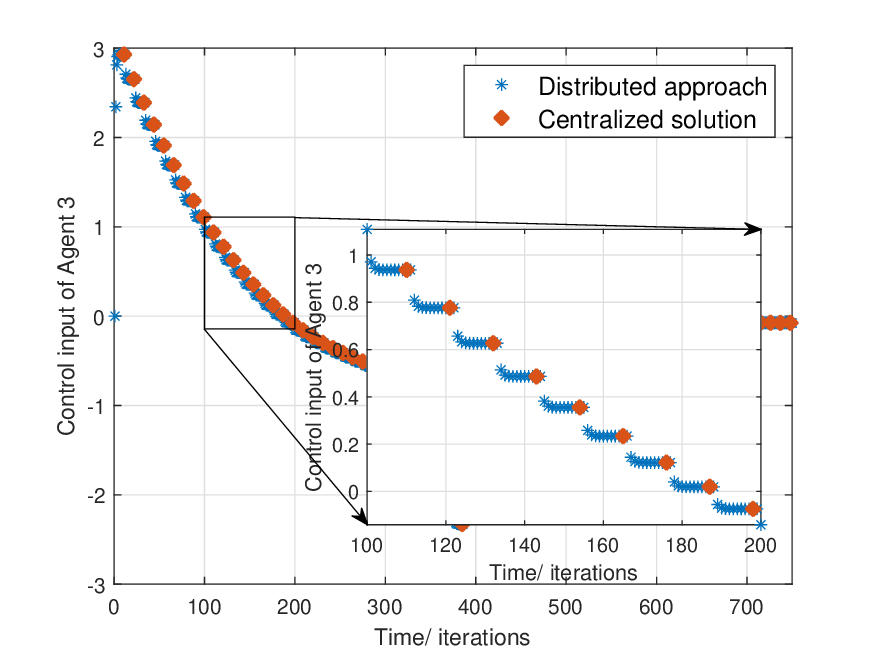}
    \caption{Control input of agent 3; centralised and distributed solution with 10 iterations/ sampling interval.}
        \label{continput3}
\end{minipage}      
\end{figure}

% \begin{figure}[!ht]
%       \centering
%       \framebox{\parbox{3in}{
%   \includegraphics[width=0.43\textwidth]{figures/uagent1.eps}}}
%         \caption{Control input of agent one; centralised solution and distributed solution with 10 iterations per sampling interval.}
%         \label{continput1}
% \end{figure}
% \begin{figure}[!ht]
%       \centering
%       \framebox{\parbox{3in}{
%   \includegraphics[width=0.43\textwidth]{figures/uagent3.eps}}}
%         \caption{Control input of agent three; centralised solution and distributed solution with 10 iterations per sampling interval.}
%         \label{continput3}
% \end{figure}

The last two plots compare control input progression achieved through the centralized solution in \eqref{udagger} and the distributed approach in \eqref{gs} with 10 iterations per sampling interval. Each plot displays the $x$-direction evolution, with blue stars indicating the distributed solution and orange diamonds representing the centralized solution. In Fig. \ref{continput3}, the distributed scheme quickly converges to the centralized solution, despite a small initial gap. Meanwhile, Fig. \ref{continput1} reveals that agent one's control input convergence is slower due to limited interaction with only two neighbors. All visualizations were generated using the code in \cite{Prima21}.

% Move this section on the full version
\subsection{Cost Comparison}

Furthermore, it is interesting to compare the costs between the original problem solved on the nodes that result in the Nash equilibrium and the reformulated problem directly solved on the edge of a network graph. Although the first problem cannot be implemented in a distributed way, we simulated it solely for cost comparison purposes, as follows. After obtaining the solution $u^{i\star}_k$ for the Nash equilibrium in \eqref{unash} and its corresponding closed-loop trajectory $x^\star_{k+1}$ in \eqref{closedloopNash}, we can derive the relative state $\bar{z}_k$ and the relative inputs $\bar{a}^m_k$ from the coupled problem, as shown in \eqref{newstate}. These values are then used to compute the cost in \eqref{cost3} and are denoted by $\bar{J}_{Nash}^m$. In contrast, for the decoupled problem, the values of $z^m_k$ and $a^m_k$ are used to compute the cost in \eqref{cost3}, resulting in $\bar{J}^m$.

\begin{table}[!h]
    \centering
    \begin{tabular}{c|c|c|c|c|c||c}
        & Edge 1 & Edge 2 & Edge 3 & Edge 4 & Edge 5 & Total\\
        \hline \hline
        $\bar{J}^m_{Nash}$ & 481 & 389 & 723 & 72 & 591 & 2256\\
        \hline
        $\Bar{J}^m$ & 293 & 127 & 667 & 13 & 554 & 1654\\
        \hline \hline
    \end{tabular}
    \caption{Cost Comparison for Each Edge Agent Using Solutions Obtained From Coupled and Decoupled Problem}
    \label{tab1}
\end{table}

Table \ref{tab1} shows that the global cost achieved by directly running the problem on the edge, $\sum_m^M \bar{J}^m = 1654$, is significantly lower than the total cost obtained by running the problem on the nodes and then transforming it, i.e., $\sum_m^M \bar{J}^m_{Nash} = 2256$. As expected, the solution attained on the edge system is Pareto optimal.

\section{CONCLUSIONS}
This article addresses the challenge of guiding a group of $N$ agents from their initial position to a desired formation while avoiding collisions with neighboring agents. The original problem is formulated as an LQDTG with a coupled SDRDE, which cannot be solved  in a distributed fashion. To address this issue, a distributed approach is proposed. This approach is based on a fictitious MAS that operates on the edges of the graph rather than the nodes. The technique incorporates relative soft constraints on the edges to prevent collisions and requires the solution of a decoupled SDRDE, using a receding horizon technique. The proposed method leverages a distributed steepest descent algorithm to map the relative control inputs to the actual physical control inputs, resulting in a simple vector-matrix multiplication per iteration, in contrast to an iterative approach often used in distributed MPC, which requires solving an optimization problem in each sampling interval. The efficacy of the proposed method is demonstrated through simulation results.

% \addtolength{\textheight}{-12cm}   % This command serves to balance the column lengths
                                  % on the last page of the document manually. It shortens
                                  % the textheight of the last page by a suitable amount.
                                  % This command does not take effect until the next page
                                  % so it should come on the page before the last. Make
                                  % sure that you do not shorten the textheight too much.

%%%%%%%%%%%%%%%%%%%%%%%%%%%%%%%%%%%%%%%%%%%%%%%%%%%%%%%%%%%%%%%%%%%%%%%%%%%%%%%%

%%%%%%%%%%%%%%%%%%%%%%%%%%%%%%%%%%%%%%%%%%%%%%%%%%%%%%%%%%%%%%%%%%%%%%%%%%%%%%%%

%%%%%%%%%%%%%%%%%%%%%%%%%%%%%%%%%%%%%%%%%%%%%%%%%%%%%%%%%%%%%%%%%%%%%%%%%%%%%%%%
\section*{APPENDIX}

\subsection*{Proof of Theorem 1}
This proof extends the methodology presented in \cite{Tamer} to the case of a state-dependent weighting matrix. Referring back to the dynamics \eqref{dyn1} and the cost \eqref{cost1}, assuming $Q^i_T(x_T) = Q^i(x_T)$, we can now rewrite it as follows:
\begin{align*}
    J^i = \frac{1}{2} \sum_{k=0}^{T} \Big( x_k^TQ^i(x_k)x_k + \sum_{j=1}^N u^{i^T}_kR^{ij}u^i_k\Big) 
\end{align*}
To maintain consistency with the indexing used in \cite{Tamer}, we will index the state costs related to the next time step.
\begin{align}\label{costmod}
    J^i = \frac{1}{2}  \sum_{k=0}^{T} \Big( x_{k+1}^TQ^i(x_{k+1})x_{k+1} + \sum_{j=1}^N u^{i^T}_kR^{ij}u^i_k\Big) 
\end{align}
To begin, it is important to observe that the Hamiltonian function for the provided equations \eqref{dyn1} and \eqref{costmod} is as follows:
\begin{align}\label{ham}
    H^i_k &= \frac{1}{2}  \Big( x_{k+1}^TQ^i(x_{k+1})x_{k+1} + \sum_{j=1}^N u^{i^T}_kR^{ij}u^i_k\Big) \nonumber\\
    &\hspace{2pc}+ \lambda_{k+1}^{i^T}\Big( Fx_k + \sum_{j=1}^N G^j u^j_k\Big)
\end{align}
Considering that $Q^i(x_{k+1})\geq 0$ and $R^{ii}>0$, setting the derivative (in control) of \eqref{ham} to zero and assuming that states and costates are fixed to their optimal values lead to the following 
\begin{align}\label{ustar}
 u^{i\star}_k &= -R^{{ii}^{-1}} G^{i^T}\Big( \lambda^i_{k+1} + Q^i(x_{k+1})x^\star_{k+1}+ (I_N \otimes x_{k+1}^T) \nonumber\\
 &\begin{bmatrix}x^{\star^T}_{k+1}
		\frac{\partial Q^i(x_{k+1})}{\partial x^1_{k+1}}x^{\star}_{k+1},&...,&x^{\star^T}_{k+1}\frac{\partial Q^i(x_{k+1})}{\partial x^N_{k+1}}x^{\star}_{k+1}
		\end{bmatrix}^T \Big).
\end{align}
Moreover, the equation for the costate (difference) is as follows
\begin{align} \label{costate}
    \lambda^i_k &= F^T\Big( \lambda^i_{k+1} +  Q^i(x_{k+1})x^\star_{k+1}+ (I_N \otimes x_{k+1}^T) \nonumber\\
 &\begin{bmatrix}x^{\star^T}_{k+1}
		\frac{\partial Q^i(x_{k+1})}{\partial x^1_{k+1}}x^{\star}_{k+1},&...,&x^{\star^T}_{k+1}\frac{\partial Q^i(x_{k+1})}{\partial x^N_{k+1}}x^{\star}_{k+1}
		\end{bmatrix}^T \Big),
\end{align}
with $\lambda^i_{T+1} = 0$, where
\begin{align}\label{cstate}
    x^\star_{k+1} = Fx_k^\star + \sum_{j=1}^N G^j u^{j\star}_k, \hspace{1pc} x_1^\star = x_1.
\end{align}
Starting with $k=T$, the expression for \eqref{ustar} simplifies to: 
\begin{align}\label{ustar2}
     u^{i\star}_k = -R^{{ii}^{-1}} G^{i^T}P^i(x_{T+1})x^\star_{T+1},
\end{align}
by first premultiplying both sides with $G^i$ and then summing over $i\in N$, and also utilizing \eqref{cstate} and \eqref{lamb}, we arrive at the following 
\begin{align*}
    x^\star_{T+1} - Fx^\star_T = (I-\Lambda_T)x^\star_{T+1}
\end{align*}
which further yields the unique relation
\begin{align*}
    x^\star_{T+1} = \Lambda_T^{-1}Fx^\star_T
\end{align*}
which is precisely \eqref{closedloopNash} for $k=T$. After substituting this relation into \eqref{ustar2}, we obtain \eqref{unash} for $k=T$.

Let's proceed to prove by induction that the unique solution set of \eqref{ustar}-\eqref{cstate} is given by \eqref{unash}-\eqref{closedloopNash} and 
\begin{align}\label{lamb2}
    \lambda^{i\star}_k = F^TP^i_{k+1}(x_{k+1})x^\star_{k+1}.
\end{align}
To achieve that, we will proceed step by step from the final time towards the initial time. For the inductive step, we start with \eqref{ustar} and substitute \eqref{lamb2}, followed by \eqref{closedloopNash}.
\begin{align*}
     u^{i\star}_k &= -R^{{ii}^{-1}} G^{i^T}\Big( \lambda^{i\star}_{k+1} + Q^i(x_{k+1})x^\star_{k+1}+ (I_N \otimes x_{k+1}^T) \nonumber\\
                     &\hspace{1pc}\begin{bmatrix}x^{\star^T}_{k+1}
                		\frac{\partial Q^i(x_{k+1})}{\partial x^1_{k+1}}x^{\star}_{k+1}&...&x^{\star^T}_{k+1}\frac{\partial Q^i(x_{k+1})}{\partial x^N_{k+1}}x^{\star}_{k+1}
                		\end{bmatrix}^T \Big)\\
	  &= -R^{{ii}^{-1}} G^{i^T}\Big(  F^TP^i_{k+2}(x_{k+2})x^\star_{k+2}+ Q^i(x_{k+1})x^\star_{k+1}\nonumber\\
	  &\hspace{1pc} +(I_N \otimes x_{k+1}^T)\nonumber\\
	  &\hspace{1.1pc}\begin{bmatrix}x^{\star^T}_{k+1}
                		\frac{\partial Q^i(x_{k+1})}{\partial x^1_{k+1}}x^{\star}_{k+1}&...&x^{\star^T}_{k+1}\frac{\partial Q^i(x_{k+1})}{\partial x^N_{k+1}}x^{\star}_{k+1}
                		\end{bmatrix}^T \Big)\\
                		&= -R^{{ii}^{-1}} G^{i^T}\Big(  F^TP^i_{k+2}(x_{k+2})\Lambda^{-1}_{k+1}Fx^\star_{k+1}\nonumber\\
	  &\hspace{1pc}+ Q^i(x_{k+1})x^\star_{k+1} +(I_N \otimes x_{k+1}^T)\nonumber\\
	  &\hspace{1.1pc}\begin{bmatrix}x^{\star^T}_{k+1}
                		\frac{\partial Q^i(x_{k+1})}{\partial x^1_{k+1}}x^{\star}_{k+1}&...&x^{\star^T}_{k+1}\frac{\partial Q^i(x_{k+1})}{\partial x^N_{k+1}}x^{\star}_{k+1}
                		\end{bmatrix}^T \Big)\\
      &= -R^{{ii}^{-1}} G^{i^T}P^i_{k+1}(x_{k+1})x^\star_{k+1},
\end{align*}
this process results in the same expression as \eqref{unash}. Next, we examine \eqref{cstate} and substitute the previous result obtained.
\begin{align*}
     x^\star_{k+1} &= Fx_k + \sum_{j=1}^N G^j u^{j\star}_k\\
     &= Fx_k - \sum_{j=1}^N G^jR^{{ii}^{-1}} G^{i^T}P^i_{k+1}(x_{k+1})x^\star_{k+1}\\
     \Big[I &+ \sum_{j=1}^N G^jR^{{ii}^{-1}} G^{i^T}P^i_{k+1}(x_{k+1})\Big]x^\star_{k+1}= Fx_k\\
     x^\star_{k+1} &= \Lambda_k^{-1}Fx_k.
\end{align*}
We have now confirmed the validity of \eqref{unash} and \eqref{closedloopNash}. To derive \eqref{lamb2}, we assume it holds up to $k+1$ and verify its applicability by substitution into \eqref{costate} for time step $k$.
\begin{align*} 
    \lambda^{i\star}_k &= F^T\Big( \lambda^{i\star}_{k+1} +  Q^i(x_{k+1})x^\star_{k+1}+ (I_N \otimes x_{k+1}^T) \nonumber\\
 &\hspace{1pc}\begin{bmatrix}x^{\star^T}_{k+1}
		\frac{\partial Q^i(x_{k+1})}{\partial x^1_{k+1}}x^{\star}_{k+1},&...,&x^{\star^T}_{k+1}\frac{\partial Q^i(x_{k+1})}{\partial x^N_{k+1}}x^{\star}_{k+1}
		\end{bmatrix}^T \Big)\\
		&= F^T\Big( F^TP^i_{k+2}(x_{k+2})x^\star_{k+2} +  Q^i(x_{k+1})x^\star_{k+1}+ \nonumber\\
		&\hspace{1pc}(I_N \otimes x_{k+1}^T) \nonumber\\
 &\hspace{0.8pc}\begin{bmatrix}x^{\star^T}_{k+1}
		\frac{\partial Q^i(x_{k+1})}{\partial x^1_{k+1}}x^{\star}_{k+1},&...,&x^{\star^T}_{k+1}\frac{\partial Q^i(x_{k+1})}{\partial x^N_{k+1}}x^{\star}_{k+1}
		\end{bmatrix}^T \Big)\\
		&= F^T\Big( F^TP^i_{k+2}(x_{k+2})\Lambda^{-1}_{k+1}Fx^\star_{k+1} +  Q^i(x_{k+1})x^\star_{k+1}+ \nonumber\\
		&\hspace{1pc}(I_N \otimes x_{k+1}^T) \nonumber\\
 &\hspace{0.8pc}\begin{bmatrix}x^{\star^T}_{k+1}
		\frac{\partial Q^i(x_{k+1})}{\partial x^1_{k+1}}x^{\star}_{k+1},&...,&x^{\star^T}_{k+1}\frac{\partial Q^i(x_{k+1})}{\partial x^N_{k+1}}x^{\star}_{k+1}
		\end{bmatrix}^T \Big)\\
		&= F^TP^i_{k+1}(x_{k+1})x^\star_{k+1},
\end{align*}
which agrees with \eqref{lamb2}, thereby completing the induction process.

% \section*{ACKNOWLEDGMENT}

%%%%%%%%%%%%%%%%%%%%%%%%%%%%%%%%%%%%%%%%%%%%%%%%%%%%%%%%%%%%%%%%%%%%%%%%%%%%%%%%

\bibliographystyle{IEEEtran}
\bibliography{IEEEexample}

\end{document}